\documentclass{aims}
\usepackage{amsmath}
 \usepackage[colorlinks=true]{hyperref}
\hypersetup{urlcolor=blue, citecolor=red}
\usepackage{amssymb,amstext}
\usepackage{capt-of}
\usepackage{textcomp}

\def\currentvolume{X}
 \def\currentissue{X}
  \def\currentyear{200X}
   \def\currentmonth{XX}
    \def\ppages{X--XX}
     \def\DOI{10.3934/xx.xx.xx.xx}

\newtheorem{theorem}{Theorem}[section]
\newtheorem{corollary}{Corollary}

\newtheorem{lemma}[theorem]{Lemma}

\theoremstyle{definition}
\newtheorem{definition}[theorem]{Definition}
\newtheorem{remark}{Remark}

\newcommand{\ep}{\varepsilon}

\def \proof{{\noindent {\bf Proof.} }}
\def \endproof{{\hfill $\Box$ \medskip}}

\newcommand{\sgn}{\mbox{\rm sgn}}

\def \N{\mathbb{N}}
\def \R{\mathbb{R}}

\def \I{\mathrm{I}}
\def \II{\mathrm{I\hspace{-.1em}I}}
\def \III{\mathrm{I\hspace{-.1em}I\hspace{-.1em}I}}
\def \IV{\mathrm{I\hspace{-.1em}V}}
\def \V{\mathrm{V}}

\def \la{\left\langle}
\def \ra{\right\rangle}
\def \ld{\left\|}
\def \rd{\right\|}
\def \lr{\left(}
\def \rr{\right)}

\title[Fast reaction limit] 
      {Fast reaction limit of reaction-diffusion systems}

\author[Murakawa]{}

\subjclass{Primary: 35K57, 35K65; Secondary: 35K51, 35R35, 80A22, 80A30.}
 \keywords{Reaction-diffusion system, fast reaction limit, singular limit, nonlinear diffusion problem, degenerate parabolic equation, $L^p$ convergence.}

 \email{murakawa@math.kyushu-u.ac.jp, hideki.murakawa@gmail.com}

\thanks{This work was supported by JSPS KAKENHI Grant nos. 26287025, 15H03635 and 17K05368. Most of the work was performed during a visit of the author to Imperial College London thanks to JST CREST Grant No. JPMJCR14D3. 
The support of JST and the hospitality of Imperial College London are warmly acknowledged.
}

\begin{document}
\maketitle

\centerline{\scshape Hideki Murakawa$^*$}
\medskip
{\footnotesize
 \centerline{Faculty of Mathematics, Kyushu University}
   \centerline{744 Motooka, Nishiku, Fukuoka 819-0395, Japan}
}
\medskip

\bigskip


\begin{abstract}
Singular limit problems of reaction-diffusion systems have been studied in cases where the effects of the reaction terms are very large compared with those of the other terms. Such problems appear in literature in various fields such as chemistry, ecology, biology, geology and approximation theory. 
In this paper, we deal with the singular limit of a general reaction-diffusion system including many problems in the literature. 
We formulate the problem, derive the limit equation and establish a rigorous mathematical theory. 
\end{abstract}

\section{Introduction}

Reaction-diffusion systems have long attracted a great deal of attention due to their high expressiveness of spatio-temporal behaviours appearing in chemistry and biology, and the richness of the structure of their solutions. 
In some problems, systems include reaction terms which are very fast comparing with the other terms. 
Singular limit of such systems as the reaction rates become extremely large, that is called a {\it fast reaction limit} or a {\it instantaneous reaction limit}, has been extensively studied in many fields of applications, e.g., 
diffusive irreversible chemical reactions~\cite{evans}, 
diffusive reversible chemical reactions~\cite{bh,ehmo}, 
spatial segregation of competing species~\cite{dhmp,ehhp}, 
invasion of bacteria in burn wounds~\cite{hkr}, 
reactive solute transport in porous media~\cite{knabner}, 
precipitation and dissolution reaction related to geological disposal of radioactive waste~\cite{behhh}, 
Stefan problem with phase relaxation~\cite{visintin_general}, 
reaction-diffusion system approximation to degenerate parabolic equations~\cite{murakawa3}, 
etc. 
Free boundaries appear in the fast reaction limit of all these problems. 
By analyzing the fast reaction limit and the ensuing free boundaries, we can understand, e.g., how biological species segregate their habitats, how the barriers that provide safety for geological waste disposal are eroded.

The reaction-diffusion systems in all of the above references can be summarized in the following system: 
\begin{equation}
\label{rd0}
\left\{
\begin{aligned}
&\frac{\partial u}{\partial t} = d_1 \Delta u +f_1(u,v) -kF(u,v)& & \mbox{in } Q_T:=\Omega \times (0,T),  \\
&\frac{\partial v}{\partial t} = d_2 \Delta v +f_2(u,v) +c kF(u,v)& & \mbox{in } Q_T,  
\end{aligned}
\right.
\end{equation}
where $d_1$, $d_2$ are constants satisfying $d_1>0$, $d_2\ge 0$, and $f_i$ $(i=1,2)$ and $F$ are nonlinear functions, a constant $c$ can be positive or negative, $k$ is a large positive constant, $\Omega\subset \mathbb{R}^d\ (d\in \mathbb{N} )$ is a bounded domain with smooth boundary $\partial \Omega$, $T$ is a positive constant. 
The unknown functions $u=u(x,t)$, $v=v(x,t)$ stand for concentrations of chemicals in the cases of chemical reactions and for population densities of species in the case of population ecology at position $x\in \Omega$ and time $t\in [0,T)$. 

Let us introduce a typical example briefly. Evans~\cite{evans} considered the following system of diffusive irreversible chemical reactions\footnote{Here, we set $k:=1/\varepsilon$ and $v:=-v^\varepsilon$, where $\varepsilon$ and $v^\varepsilon$ are notations used in~\cite{evans}.}: 
\begin{equation}
\label{irreversible}
\left\{
\begin{aligned}
&\frac{\partial u}{\partial t} = d_1 \Delta u -kuv & & \mbox{in } Q_T,  \\
&\frac{\partial v}{\partial t} = d_2 \Delta v  -kuv & & \mbox{in } Q_T
\end{aligned}
\right.
\end{equation}
with non-negative initial data. 
Here, $u$ and $v$ denote the concentrations of the two reacting chemicals which combine to form another product. 
We are interested in the behaviours of the chemicals when the reaction rate $k$ is very large. 
Let $(u^k,v^k)$ be the solution of \eqref{irreversible}. 
It has been shown that 
the limit functions of $u^k$ and $v^k$ as $k$ tends to infinity are given by positive $z^+:=\max\{ z,0\}$ and negative $z^-:=\max\{ -z,0\}$ parts of the solution $z$ of the following nonlinear diffusion equation, respectively: 
\begin{equation}
\label{limit0}
\frac{\partial z}{\partial t} = \Delta \beta(z),
\end{equation}
where $\beta(z):=d_1\max\{z,0\}+d_2\min\{z,0\}$. 
This equation~\eqref{limit0} is a weak form of a free boundary problem so-called the Stefan problem. 
The behaviours of the free boundary is well-known and one can estimate the errors between $(u^k,v^k)$ and $(z^+,z^-)$~\cite{murakawa2}. Therefore, we can understand the behaviours of solutions near the strongly reacting front. 
This sort of fast reaction limit, especially, with the logistic type of $f_1$ and $f_2$, has been widely studied (see, e.g., \cite{dhmp,ehhp}).

Like this example, if $c$ is negative, initial data are non-negative and $F$ is non-decreasing in both variables with $F(a,b)=0$ iff $ab=0$, then \eqref{rd0} is called irreversible reaction-diffusion system. 
If $c$ is positive and $F(u,v)$ is non-decreasing in $u$ and non-increasing in $v$, then \eqref{rd0} is called reversible reaction-diffusion system. 
The irreversible and the reversible systems have been handled separately because it has been thought that different techniques are required. 
However, the irreversible system can be written as the reversible system formally by replacing $v$ with $-v$. 
Precisely speaking, for $F(u,v)=uv$ as in \eqref{irreversible}, $F(u,-v)$ does not satisfy the monotonicity assumption. 
If we replace $F(u,v)=uv$ with 
\begin{equation*}
F(u,v) =
\left\{
\begin{aligned}
& uv & & \mbox{if } u\ge 0,\ v\ge 0,  \\
&u & & \mbox{if } u<0,\ u\le v, \\
&v & & \mbox{if } u<0,\ u> v, 
\end{aligned}
\right.
\end{equation*}
then $F(u,-v)$ satisfies the monotonicity assumption and this does not change the dynamics within the invariant set $\{ u,v\ge 0\}$. 
Thus, the reaction-diffusion systems in all of the above references can be reformulated as the reversible system.

In the present paper, we consider the fast reaction limit of the following general system: 
\begin{equation}
\label{rd}
\left\{
\begin{aligned}
&\frac{\partial u}{\partial t} = d_1 \Delta u +f_1(u,v) -kF(u,v) & & \mbox{in } Q_T,  \\
&\frac{\partial v}{\partial t} = d_2 \Delta v +f_2(u,v) +kF(u,v) & & \mbox{in } Q_T,  \\
&\frac{\partial u}{\partial \nu}=\frac{\partial v}{\partial \nu}= 0 & & \mbox{on } \partial \Omega \times (0,T), \\
&u(0) = u_0, \quad v(0)=v_0 & & \mbox{in } \Omega.
\end{aligned}
\right.
\end{equation}
Here, $\nu$ is the outward normal unit vector to $\partial \Omega$ and $u_0$, $v_0$ are given initial data. 
We derive the limit problem as $k$ tends to infinity and prove the convergence under a broad framework including all of the above references. 

Fast reaction limit of three- or multi-component systems have been studied in \cite{himn,hm,in,m_lscds,mn}. The structures of these problems are similar to our framework essentially. 
In our setting~\eqref{rd}, 
the fast reaction terms in each component $u$ and $v$ are represented by the same function $F$. Fast reaction limits of systems with different fast reaction terms in equations for $u$ and $v$ are still far from being well understood. 
There are some researches on such problems, see, e.g., \cite{ctv1,immh}. 
Quite recently, Iida, Ninomiya and Yamamoto reviewed the recent development of fast reaction limits and reaction-diffusion system approximations~\cite{iny}.

The organization of this paper is as follows. 
In the next section, we state notations and definitions, and give our main results. 
In Section~3, a priori estimates are provided in $L^2$ framework. 
The main result in $L^p$ ($p\in [1,2)$) space is proved in Section~4. The $L^\infty$ estimate and the main result in $L^p$ ($p\ge 1$) space are given in Section~5. 
Concluding remarks are made in the final section of the paper.

\section{Assumptions and main results}
\label{Sec_results}

In this paper, we deal with the fast reaction limit of \eqref{rd} under the following general setting. 
\renewcommand{\theenumi}{H\arabic{enumi}}
\renewcommand{\labelenumi}{(\theenumi)}
\begin{enumerate}
\item \label{assum_d}
$d_1>0$, $d_2 \ge 0$. 
\item \label{assum_f12}
For $i=1,2$, $f_i:\R^2 \to \R$ is Lipschitz continuous. 
\item \label{assum_F}
$F:\R^2 \to \R$ is Lipschitz continuous. \\
$F(u,v)$ is non-decreasing in $u$ and non-increasing in $v$. \\
There exists a maximal monotone $\alpha : \R  \to 2^\R$ such that 
\begin{equation}
\label{rel_F_alpha}
F(u,v)=0 \Longleftrightarrow v\in \alpha(u) \qquad \mbox{for all } u,v\in \R, 
\end{equation}
$\alpha$ grows at least and at most linearly at infinity, more precisely, there exist positive constants $C_i$ $(i=1,2,3,4)$ such that for all $u\in \R$ and $v\in \alpha (u)$, the following relation holds.
\begin{equation}
\label{propertyofalpha1}
C_1|u|-C_2 \le |v| \le C_3|u|+C_4.
\end{equation}
\end{enumerate}
Since $\alpha$ is maximal monotone, \eqref{rel_F_alpha} yields 
\begin{align*}
F(u,v)>0 &\Longleftrightarrow v< \mbox{inf}\, \alpha(u), \\
F(u,v)<0 &\Longleftrightarrow v> \mbox{sup}\, \alpha(u)
\end{align*}
for all $u,v\in \R$. 
Condition~\eqref{propertyofalpha1} implies 
\begin{equation*}
\frac{1}{C_3}|v|-\frac{C_4}{C_3} \le |u| \le \frac{1}{C_1}|v|+\frac{C_2}{C_1}
\end{equation*}
for all $v\in \R$ and $u\in \alpha^{-1} (v)$. 

Let us formally derive the limit functions $(u^*,v^*)$ of the solution $(u^k,v^k)$ of \eqref{rd} as $k$ goes to infinity. 
It follows from the equation for $u^k$ that 
\[
F(u^k,v^k)=\frac{1}{k}\left(
-\frac{\partial u^k}{\partial t} + d_1 \Delta u^k +f_1(u^k,v^k)
\right).
\]
If the term in the right parenthesis is bounded with respect to $k$ (or lower order than the first order) in a sense, passing to the limit in $k$ yields $F(u^*,v^*)=0$, and thus, $v^*\in \alpha(u^*)$. 
Setting $z = u^*+v^*$, we have $u^*=(I+\alpha)^{-1}z$, $v^*=(I-(I+\alpha)^{-1})z$ and $d_1u^*+d_2v^*=(d_2I+(d_1-d_2)(I+\alpha)^{-1})z$. 
Here and hereafter, $I$ denotes the identity function. 
Since $\alpha$ is maximal monotone, $(I+\alpha)^{-1}$ is a non-decreasing contraction mapping. 
On the other hand, we deduce from the equations for $u^k$ and $v^k$ that 
\[
\frac{\partial }{\partial t}\lr u^k+v^k\rr = \Delta \lr d_1u^k +d_2v^k\rr +f_1(u^k,v^k)+f_2(u^k,v^k). 
\]
Letting to the limit in $k$ and substituting $u^*$ and $v^*$, we have an equation for $z$. 
In summary, we obtain the following assertion.

The limit functions $u^*, v^*$ of the solutions of \eqref{rd} as $k$ tends to infinity can be represented by $u^*=(I+\alpha)^{-1}z$ and $v^*=(I-(I+\alpha)^{-1})z$ using a (weak) solution $z$ of the following nonlinear diffusion problem. 
\begin{equation}
\label{limiteq}
\left\{
\begin{aligned}
&\frac{\partial z}{\partial t} = \Delta \beta(z)+f(z) & & \mbox{in } Q_T,  \\
&\frac{\partial \beta(z)}{\partial \nu} = 0 & & \mbox{on } \partial \Omega \times (0,T), \\
&z(0) = u_0+v_0 & & \mbox{in } \Omega.
\end{aligned}
\right.
\end{equation}
Here, $\beta$ and $f$ are defined as 
\begin{equation}
\label{def_beta_f}
\begin{aligned}
\beta(s) & = (d_2I+(d_1-d_2)(I+\alpha)^{-1})s, \\
f(s) &= (f_1+f_2)((I+\alpha)^{-1}s, (I-(I+\alpha)^{-1})s)
\end{aligned}
\qquad \mbox{for } s\in \R. 
\end{equation}
We note that $\beta$ and $f$ are Lipschitz functions, and $\beta$ is non-decreasing because $(I+\alpha)^{-1}$ is contraction. 
If $d_2=0$, then \eqref{limiteq} can be degenerate, and \eqref{limiteq} involves the Stefan or the porous medium type free boundary problems.

Our assertion is that the fast reaction limit of \eqref{rd} is described by the nonlinear diffusion problem \eqref{limiteq}. Conversely, the solution of the nonlinear diffusion problems of the type \eqref{limiteq} can be approximated by a solution of a reaction-diffusion system. 
In this context, we consider a general initial datum $(u_0, v_0)$ for \eqref{limiteq} and its approximation $(u_0^k, v_0^k)$ as an initial datum of \eqref{rd}. 
We impose the following assumptions on the initial data: 

\renewcommand{\theenumi}{H\arabic{enumi}}
\renewcommand{\labelenumi}{(\theenumi)}
\begin{enumerate}
\setcounter{enumi}{3}
\item \label{assum_init}
\begin{align*}
&u_0 \in L^2(\Omega),\ v_0 \in L^2(\Omega),\ \Delta u_0 \in C_0^0(\overline{\Omega})^*,\ 
d_2 \Delta v_0 \in C_0^0(\overline{\Omega})^*, \\
&F(u_0,v_0)=0 \quad \mbox{a.e. in } \Omega. 
\end{align*}
For all $k>0$, 
\begin{align}
&u_0^k \in H^1(\Omega),\ v_0^k \in H^1(\Omega),\ \Delta u_0^k \in L^1({\Omega}),\ 
d_2 \Delta v_0^k \in L^1({\Omega}), \notag \\
&\frac{\partial u_0^k}{\partial \nu}=\frac{\partial v_0^k}{\partial \nu} = 0 \quad \mbox{on } \partial \Omega, \notag \\
&F(u_0^k,v_0^k)\le C_5/k \quad \mbox{a.e. in } \Omega, \label{assum_init_C5}\\
&\ld u_0^k \rd_{L^2(\Omega)} + \ld v_0^k \rd_{L^2(\Omega)} + \ld \Delta u_0^k\rd_{L^1(\Omega)}+d_2 \ld \Delta u_0^k\rd_{L^1(\Omega)} \le C_6, \notag \\
&u_0^k \to u_0, \ v_0^k \to v_0 \mbox{ in } L^2(\Omega) \mbox{ as } k\to \infty. \notag 
\end{align}
Here $C_5$ and $C_6$ are positive constant independent of $k$. 
\end{enumerate}

Under these assumptions, \eqref{rd} has the unique solution $(u,v)$ such that 
${u}\in H^1(0,T;L^2(\Omega)) \cap L^\infty(0,T;H^1(\Omega))$ and $\Delta {u} \in L^2(Q_T)$, and $v$ has the same regularities if $d_2>0$, $v \in W^{1,\infty}(0,T;L^2(\Omega)) \cap L^\infty(0,T;H^1(\Omega))$ if $d_2=0$.

We denote by $\langle \cdot,\cdot\rangle $ the duality pairing between $H^{1}(\Omega)^*$ and $H^1(\Omega)$ and by $\lr \cdot, \cdot \rr$ the inner product in $L^2(\Omega)$.

Problem \eqref{limiteq} is understood in the following weak sense. 
\begin{definition}
\label{weak_limiteq}
A function ${z}\in L^\infty(0,T;L^2(\Omega))\cap H^1(0,T;H^{1}(\Omega)^*)$ such that ${\beta}({z}) \in L^2(0,T;H^1(\Omega))$ is said to be a weak solution of \eqref{limiteq} if it fulfils 
\begin{equation*}
-\int_0^T\left\langle \frac{\partial \varphi}{\partial t}, z \right\rangle 
+\int_0^T \lr \nabla \beta({z}),\nabla \varphi\rr
 = \lr u_0+v_0, \varphi(0) \rr + \int_0^T\lr f({z}), \varphi\rr
\end{equation*}
for all functions $\varphi \in L^2(0,T;H^1(\Omega))$ such that $\varphi (T)=0$. 
\end{definition}

Existence of the weak solution of the nonlinear diffusion equations of type \eqref{limiteq} is well-known (see, e.g., \cite{crank,verdi2}). However, uniqueness of the weak solution is not known in general. The weak solution is unique if one of the following additional assumptions holds (see, e.g., \cite{acp,jerome}):
\begin{itemize}
\item $d_2>0$, 
\item 
$\beta$ is strictly monotone, 
\item there is a Lipschitz function $\tilde{f}$ such that $f(z)=\tilde{f}(\beta(z))$.
\end{itemize}

The main result in this paper is as follows. 
\begin{theorem}
\label{Thm_conv}
Assume that \eqref{assum_d}--\eqref{assum_init} are satisfied. 
Let $\beta$ and $f$ be the functions defined in \eqref{def_beta_f}. 
Let $(u^k,v^k)$ be the weak solutions of \eqref{rd} with an initial datum $(u_0^k,v_0^k)$. 
Then, there exist subsequences $\{u^{k_n}\}$ and $\{v^{k_n}\}$ of $\{u^k\}$ and $\{v^k\}$, respectively, and a weak solution $z$ of \eqref{limiteq} with the initial datum $u_0+v_0$ such that 
\begin{align*}
{z} \in L^\infty(0,T;L^2(\Omega)), \quad 
\beta(z) \in {L^2(0,T;H^1(\Omega))} \\
(z \in {L^2(0,T;H^1(\Omega))} \mbox{ if }  d_2>0),
\end{align*}
\[
\begin{array}{ll}
u^{k_n} \to (I+\alpha)^{-1}z & \mbox{strongly in } L^s(Q_T)\ \forall s\in [1,2), \mbox{ a.e. in } Q_T, \\
&\mbox{weakly in } L^2(0,T;H^1(\Omega)), \\
v^{k_n} \to (I-(I+\alpha)^{-1})z & \mbox{strongly in } L^s(Q_T)\ \forall s\in [1,2), \mbox{ a.e. in } Q_T, \\
&\mbox{weakly in } L^2(Q_T)\\
&(\mbox{weakly in } L^2(0,T;H^1(\Omega) \mbox{ if } d_2>0),
\end{array}
\]
as $k_n\to \infty$. 

Moreover, if there exists a positive constant $C$ independent of $k$ such that 
\begin{equation}
\label{init_linf}
\ld u_0^k \rd_{L^\infty(\Omega)}
+
\ld v_0^k \rd_{L^\infty(\Omega)}
\le C, 
\end{equation}
then $z\in L^\infty(Q_T)$ and further subsequences converge strongly in $L^p(Q_T)$ for all $p\ge 1$. 
\end{theorem}

Throughout this paper, we denote by $C$ a generic positive constant independent of the relevant parameters, $k$ here, $k$, $\tau$ and $\xi$ in Lemmas~\ref{Lem_time_translation} and \ref{Lem_space_translation}.

\begin{remark}
When \eqref{init_linf} holds, Assumptions~\eqref{assum_f12} and \eqref{assum_F} can be weakened. 
The global Lipschitz continuities of $f_i(i=1,2)$ and $F$ can be replaced by local Lipschitz continuities, and \eqref{propertyofalpha1} is not required. 
\end{remark}

\begin{remark}
If the weak solution of \eqref{limiteq} is unique as mentioned above, whole sequences converge. 
\end{remark}

\section{A priori estimates} 

This section is devoted to a priori estimates under Assumptions \eqref{assum_d}--\eqref{assum_init}.

\begin{lemma}
\label{Lem_apriori_L2H1}
There exists a positive constant $C$ independent of $k$ such that 
\begin{align*}
\| u^k \|_{L^\infty(0,T;L^2(\Omega))\cap L^2(0,T;H^1(\Omega))}
+
\| v^k \|_{L^\infty(0,T;L^2(\Omega))}
+
d_2 \| v^k \|_{L^2(0,T;H^1(\Omega))}
\le 
C. 
\end{align*}
\end{lemma}

\proof
We use the following notation for a given nondecreasing Lipschitz continuous function $g$ with a Lipschitz constant $L_g$: 
\begin{equation*}
\Phi_g (s) = \int_0^s g(r)dr \qquad \mbox{for } s\in \R.
\end{equation*}
It is easy to see that $\Phi_g$ is convex and satisfies 
\begin{equation*}
\frac{1}{2L_g}g(s)^2 \le \Phi_g (s) \le \frac{L_g}{2}s^2 \qquad \mbox{for } s\in \R.
\end{equation*}

Multiplying the equation for $u^k$ by $u^k$, integrating both sides in space and using integration by parts, we have
\begin{equation}
\label{apriori_ubyu}
\frac{1}{2}\frac{d}{dt} \ld u^k \rd_{L^2(\Omega)}^2
+d_1\ld \nabla u^k\rd_{L^2(\Omega)}^2 
+k \lr F(u^k,v^k),u^k \rr
=
\lr f_1(u^k,v^k),u^k\rr. 
\end{equation}

Let $\{ \gamma_\delta \}_{\delta>0}$ be a sequence of nondecreasing and Lipschitz continuous functions and $\gamma_0$ be a nondecreasing single-valued function such that for all $s\in \R$ 
\begin{align}
&\left| \gamma_\delta(s) \right| 
\le \frac{1}{C_1}|s| + \frac{C_2}{C_1}, \label{alphadelta_property1} \\
& \gamma_0(s) \in \alpha^{-1}(s), \notag \\
& \gamma_\delta(s) \to \gamma_0(s) \mbox{ as } \delta \to 0. \notag
\end{align}
Here, we note that $\alpha^{-1}$ is generally multi-valued, but $\gamma_\delta$ and $\gamma_0$ are single-valued, and one can take such an approximation, e.g., Yosida approximation. 
Multiply the equation for $v^k$ by $\gamma_\delta(v^k)$, and integrate both sides in space to get 
\begin{equation}
\label{apriori_vbyalphav}
\frac{d}{dt} \int_\Omega \Phi_{\gamma_\delta}(v^k) 
+d_2\lr \nabla \gamma_\delta(v^k), \nabla v^k \rr  
-k \lr F(u^k,v^k),\gamma_\delta(v^k) \rr
=
\lr f_2(u^k,v^k),\gamma_\delta(v^k) \rr. 
\end{equation}
Combining \eqref{apriori_ubyu} and \eqref{apriori_vbyalphav}, and integrating the result in time from $0$ to $t_0 \in (0,T)$, we obtain 
\begin{align}
\frac{1}{2}&\ld u^k (t_0) \rd_{L^2(\Omega)}^2
+d_1 \ld \nabla u^k\rd_{L^2(0,t_0;L^2(\Omega))}^2 
+\int_\Omega \Phi_{\gamma_\delta}(v^k(t_0))
+d_2 \int_0^{t_0} \lr \nabla \gamma_\delta(v^k), \nabla v^k \rr  \notag \\
&+k \int_0^{t_0} \lr F(u^k,v^k),u^k -\gamma_\delta(v^k) \rr \notag \\
=&
\frac{1}{2}\ld u_0^k \rd_{L^2(\Omega)}^2
+\int_\Omega \Phi_{\gamma_\delta}(v_0^k)
+\int_0^{t_0} \left\{ \lr f_1(u^k,v^k),u^k\rr + \lr f_2(u^k,v^k),\gamma_\delta(v^k) \rr \right\}. \label{apriori_tmp1}
\end{align}
Using \eqref{alphadelta_property1} and Lipschitz continuities of $f_1$ and $f_2$, 
the second and the third terms of the right hand side of \eqref{apriori_tmp1}, which are denoted by $\II$ and $\III$, respectively, can be estimated as follows. 
\begin{align*}
\II &\le \int_\Omega \lr \frac{1}{C_1}\left| v_0^k \right|^2 + \frac{C_2}{C_1}\left| v_0^k \right| \rr 
\le C\lr 1+\ld v_0^k \rd_{L^2(\Omega)}^2 \rr \le C, \\ 
\III 
&\le C\lr 1+\int_0^{t_0} \left\{ \ld u^k \rd_{L^2(\Omega)}^2 + \ld v^k \rd_{L^2(\Omega)}^2 \right\} \rr. 
\end{align*}
The third and the fourth terms of the left hand side of \eqref{apriori_tmp1} is positive. Passing to the limit as $\delta \to 0$, we realize that the fifth term of the left hand side of \eqref{apriori_tmp1} is also positive. Thus, we have 
\begin{align}
\ld u^k (t_0) \rd_{L^2(\Omega)}^2
+ \ld \nabla u^k\rd_{L^2(0,t_0;L^2(\Omega))}^2 
+k \int_0^{t_0} \lr F(u^k,v^k),u^k -\gamma_0(v^k) \rr \notag \\
\le C\lr 1 + \int_0^{t_0} \left\{ \ld u^k \rd_{L^2(\Omega)}^2 + \ld v^k \rd_{L^2(\Omega)}^2 \right\}\rr. \label{apriori_tmp2}
\end{align}
Let $\{ \alpha_\delta \}_{\delta \ge 0}$ be a sequence of nondecreasing functions that approximates to $\alpha$ such that it has similar properties to $\gamma_\delta$. 
Multiplying the equation for $u^k$ by $\alpha_\delta (u^k)$ and the equation for $v^k$ by $v^k$, adding the results and using a similar strategy to the above, we obtain
\begin{align}
\ld v^k (t_0) \rd_{L^2(\Omega)}^2
+ d_2 \ld \nabla v^k\rd_{L^2(0,t_0;L^2(\Omega))}^2 
+k \int_0^{t_0} \lr F(u^k,v^k),\alpha_0(u^k) -v^k \rr \notag \\
\le C\lr 1+ \int_0^{t_0} \left\{ \ld u^k \rd_{L^2(\Omega)}^2 + \ld v^k \rd_{L^2(\Omega)}^2 \right\}\rr. \label{apriori_tmp3}
\end{align}
The desired estimate follows from \eqref{apriori_tmp2}, \eqref{apriori_tmp3} and the Gronwall inequality. 
Moreover, we have 
\begin{align}
&0\le k \int_0^{T} \lr F(u^k,v^k),u^k -\gamma_0(v^k) \rr \le C, \label{conv_ugamma}\\
&0\le k \int_0^{T} \lr F(u^k,v^k),\alpha_0(u^k) -v^k \rr \le C. \label{conv_alphav}
\end{align} 
\endproof

Set 
\[
\sgn(s) = 
\left\{
\begin{array}{lll}
1 & \mbox{if} & s>0,\\
0 & \mbox{if} & s=0,\\
-1 & \mbox{if} & s<0. 
\end{array}
\right.
\]
Define a Lipschitz approximation of $\sgn$ by $\sgn_n(s)=\min \left\{ \max \left\{ ns,-1 \right\}, 1\right\}$ for $s\in \R$ and $n\in \N$. 
\begin{lemma}
\label{Lem_F_sgn}
For all $u_1,u_2,v_1,v_2\in \R$, the following estimate holds. 
\[
\lr F(u_1,v_1)-F(u_2,v_2)\rr 
\lr \sgn(u_1-u_2)-\sgn (v_1-v_2) \rr \ge 0. 
\]
\end{lemma}
\proof
If either $u_1>u_2$ and $v_1>v_2$ or $u_1<u_2$ and $v_1<v_2$ hold, then $\sgn(u_1-u_2)=\sgn (v_1-v_2)$. Thus, we have the desired estimate. If $u_1\ge u_2$ and $v_1 \le v_2$, then $\sgn(u_1-u_2)\ge \sgn (v_1-v_2)$ and $F(u_1,v_1)\ge F(u_2,v_2)$ because of the property of $F$. Therefore, the desired estimate holds. In case where $u_1\le u_2$ and $v_1 \ge v_2$, the result is obtained similarly. 
\endproof

\begin{lemma}
\label{Lem_time_translation}
For all $s\in [1,2)$, there is a positive constant $C$ independent of $k$ and $\tau$ such that 
\begin{align*}
&\lr \int_0^{T-\tau}\!\!\! \int_{\Omega} \left|u^k(x,t+\tau)-u^k (x,t)\right|^s dxdt \rr^{\frac{1}{s}} \\
&
\quad + \lr \int_0^{T-\tau}\!\!\! \int_{\Omega} \left|v^k(x,t+\tau)-v^k (x,t)\right|^s dxdt\rr^{\frac{1}{s}} \le C\tau^{\frac{2-s}{s}} 
\end{align*}
for all $\tau \in (0,T)$.  
\end{lemma}

\proof
Multiplying the equation for $u^k$ by $\sgn_n (u^k-u_0^k)$ and integrating the result in space yield 
\begin{align*}
&\lr \frac{\partial u^k}{\partial t} , \sgn_n (u^k-u_0^k) \rr
-d_1 \lr \Delta (u^k-u_0^k), \sgn_n (u^k-u_0^k) \rr\notag \\
&\qquad -d_1 \lr \Delta u_0^k, \sgn_n (u^k-u_0^k) \rr 
 + k \lr F(u^k,v^k),\sgn_n (u^k-u_0^k) \rr \notag \\
&\quad =
\lr f_1(u^k,v^k),\sgn_n (u^k-u_0^k) \rr, 
\end{align*}
that is denoted by $\I+\II+\III+\IV=\V$. 
Each term can be estimated as follows. 
\begin{align*}
\II &= \lr \nabla (u^k-u_0^k), (\sgn_n)'(u^k-u_0^k)\nabla (u^k-u_0^k) \rr \ge 0, \\
\III & \le d_1 \ld \Delta u_0^k \rd_{L^1(\Omega)}\le C, \\
\IV &\le \ld f_1(u^k,v^k) \rd_{L^1(\Omega)} \ld \sgn_n (u^k-u_0^k) \rd_{L^\infty(\Omega)}
\le C.
\end{align*}
Note that $\sgn_n (u^k-u_0^k) \to \sgn (u^k-u_0^k)$ a.e. as $n\to \infty$. 
Passing to the limit in $n$, we have 
\begin{equation*}
\lr \frac{\partial u^k}{\partial t} , \sgn (u^k-u_0^k) \rr
+ k \lr F(u^k,v^k),\sgn (u^k-u_0^k) \rr 
\le C. 
\end{equation*}
The first term of the left hand side is calculated as follows. 
\[
\lr \frac{\partial u^k}{\partial t} , \sgn (u^k-u_0^k) \rr
=
\lr \frac{\partial (u^k-u_0^k)}{\partial t} , \sgn (u^k-u_0^k) \rr
=\frac{d}{dt}\ld u^k-u_0^k \rd_{L^1(\Omega)}. 
\]
Similarly, considering the multiplication of the equation for $v^k$ by $\sgn_n(v^k-v_0^k)$ yields 
\begin{equation*}
\frac{d}{dt}\ld v^k-v_0^k \rd_{L^1(\Omega)}
- k \lr F(u^k,v^k),\sgn (v^k-v_0^k) \rr 
\le C. 
\end{equation*}
Therefore, using Lemma~\ref{Lem_F_sgn} and \eqref{assum_init_C5}, we obtain 
\begin{equation*}
\frac{d}{dt}\ld u^k-u_0^k \rd_{L^1(\Omega)}
+\frac{d}{dt}\ld v^k-v_0^k \rd_{L^1(\Omega)}
\le C. 
\end{equation*}
Integrate it in time from $0$ to $\tau \in (0,T)$ to get 
\begin{equation}
\label{tt_firstestim}
\ld u^k(\tau)-u_0^k \rd_{L^1(\Omega)}
+\ld v^k(\tau)-v_0^k \rd_{L^1(\Omega)}
\le C\tau . 
\end{equation}
For $\tau\in (0,T)$, set 
\[
\begin{aligned}
\bar{u}^k(x,t)&= u^k(x,t+\tau), \\
\bar{v}^k(x,t)&= v^k(x,t+\tau)
\end{aligned}
\quad \mbox{for } x\in \Omega,\ t\in (0,T-\tau).
\]
We deduce from the equations and regularities of $u^k$ and $v^k$ that 
\begin{align}
\frac{\partial}{\partial t}(\bar{u}^k-u^k) - d_1 \Delta (\bar{u}^k-u^k) + k\lr F(\bar{u}^k,\bar{v}^k)-F(u^k,v^k) \rr \notag \\
= \lr f_1(\bar{u}^k,\bar{v}^k)-f_1(u^k,v^k) \rr, \label{tt_baru-u} \\
\frac{\partial}{\partial t}(\bar{v}^k-v^k) - d_2 \Delta (\bar{v}^k-v^k) - k\lr F(\bar{u}^k,\bar{v}^k)-F(u^k,v^k) \rr \notag \\
= \lr f_2(\bar{u}^k,\bar{v}^k)-f_2(u^k,v^k) \rr \label{tt_barv-v}
\end{align}
a.e. in $Q_{T-\tau}$. 
Multiply \eqref{tt_baru-u} by $\sgn_n (\bar{u}^k-u^k)$ and \eqref{tt_barv-v} by $\sgn_n (\bar{v}^k-v^k)$, add the results and integrate it over $\Omega$. Then, we have 
\begin{align}
&\lr \frac{\partial (\bar{u}^k-u^k)}{\partial t} , \sgn_n (\bar{u}^k-u^k) \rr
+
\lr \frac{\partial (\bar{v}^k-v^k)}{\partial t} , \sgn_n (\bar{v}^k-v^k) \rr \notag \\
& \quad -d_1 \lr \Delta (\bar{u}^k-u^k), \sgn_n (\bar{u}^k-u^k) \rr
-d_2 \lr \Delta (\bar{v}^k-v^k), \sgn_n (\bar{v}^k-v^k) \rr \notag \notag \\
& \quad + k \lr F(\bar{u}^k,\bar{v}^k)-F(u^k,v^k) ,\sgn_n (\bar{u}^k-u^k) - \sgn_n (\bar{v}^k-v^k) \rr \notag \notag \\
&=
\lr f_1(\bar{u}^k,\bar{v}^k)-f_1(u^k,v^k) ,\sgn_n (\bar{u}^k-u^k) \rr \notag \\
& \quad +
\lr f_2(\bar{u}^k,\bar{v}^k)-f_2(u^k,v^k) ,\sgn_n (\bar{v}^k-v^k) \rr. \label{tt_3}
\end{align}
Integration by parts yields that the third and the fourth terms of the left hand side are non-negative. 
In view of the Lipschitz continuities of $f_1$ and $f_2$, and Lemma~\ref{Lem_apriori_L2H1}, the right hand side is bounded by $C\lr 1+\ld \bar{u}^k-{u}^k \rd_{L^1(\Omega)} + \ld \bar{v}^k-{v}^k \rd_{L^1(\Omega)} \rr$ for a positive constant $C$ independent of $k$ and $\tau$. Passing to the limit in $n$, the fifth term of the left hand side is non-negative due to Lemma~\ref{Lem_F_sgn}. 
Therefore, we obtain 
\begin{equation}
\frac{d}{dt}\lr \ld \bar{u}^k-u^k \rd_{L^1(\Omega)}
+\ld \bar{v}^k-v^k \rd_{L^1(\Omega)} \rr
\le C\lr \ld \bar{u}^k-u^k \rd_{L^1(\Omega)}
+\ld \bar{v}^k-v^k \rd_{L^1(\Omega)} \rr. \label{tt_4}
\end{equation}
We deduce from the Gronwall inequality and \eqref{tt_firstestim} that 
\begin{align*}
&\ld \bar{u}^k-u^k \rd_{L^1(\Omega)}(t)
+\ld \bar{v}^k-v^k \rd_{L^1(\Omega)}(t) \\
& \quad \le
C\lr \ld {u}^k(\tau)-u_0^k \rd_{L^1(\Omega)}
+\ld {v}^k(\tau)-v_0^k \rd_{L^1(\Omega)} \rr \\
& \quad \le 
C\tau. 
\end{align*}
Hence, 
\begin{align*}
&\ld \bar{u}^k-u^k \rd_{L^1(Q_{T-\tau})}
+\ld \bar{v}^k-v^k \rd_{L^1(Q_{T-\tau})} 
\le
C\tau. 
\end{align*}
The desired estimate follows from the following H\"older inequality: 
\begin{equation}
\label{holder}
\ld g \rd_{L^s} 
\le \ld g \rd_{L^1}^{\frac{2-s}{s}}
 \ld g \rd_{L^2}^{\frac{2s-2}{s}}
\end{equation}
for $s\in [1,2)$ and $g\in L^2$.
\endproof

\begin{lemma}
\label{Lem_space_translation}
For $s\in [1,2)$ and sufficiently small $r>0$, there exists a positive function $\sigma$ such that $\sigma (\xi)\to 0$ as $\xi \to 0$ and the following estimate holds. 
\begin{align*}
&\lr \int_0^T\!\!\! \int_{\Omega_r} \left|u^k(x+\xi,t)-u^k (x,t)\right|^s dxdt \rr^{\frac{1}{s}}  \le \sigma(\xi), \\ 
&\lr \int_0^T\!\!\! \int_{\Omega_r} \left|v^k(x+\xi,t)-v^k (x,t)\right|^s dxdt\rr^{\frac{1}{s}} \le \sigma(\xi) 
\end{align*}
for all $\xi \in \R^N$, $|\xi|\le 2r$. Here, $\Omega_r=\{ x\in \Omega \ :\ B(x,2r)\subset \Omega \}$ and $B(x,r)$ is the open ball of radius $r$ centered at a point $x$. 
\end{lemma}

\proof
Define $\Omega_r'=\cup_{x\in \Omega_r}B(x,r)$. Remark that $\Omega_r\subset \Omega_r' \subset \Omega$. 
Let us fix a function $\psi\in C_0^\infty (\Omega_r')$ dependent only on $\Omega$ and $r$ such that 
\[
0\le \psi \le 1 \mbox{ in } \Omega_r', \quad \psi =1 \mbox{ in } \Omega_r.
\]
For $\xi \in \R^d$ with $|\xi |\le r$, set 
\[
\begin{aligned}
\tilde{u}^k(x,t)&= u^k(x+\xi,t), \\
\tilde{v}^k(x,t)&= v^k(x+\xi,t)
\end{aligned}
\quad \mbox{for } x\in \Omega_r',\ t\in (0,T).
\]
We have \eqref{tt_baru-u} and \eqref{tt_barv-v} where $\bar{u}^k$ and $\bar{v}^k$ are replaced by $\tilde{u}^k$ and $\tilde{v}^k$, respectively. 
Multiplying the equation for $\tilde{u}^k-u^k$ by $\sgn_n (\tilde{u}^k-u^k)\psi$ and the equation for $\tilde{v}^k-v^k$ by $\sgn_n (\tilde{v}^k-v^k)\psi$, adding the results and integrate it over $\Omega_r'$, 
we get similar equality to \eqref{tt_3}. 
Note that 
\begin{align*}
-& \int_{\Omega_r} \Delta (\tilde{u}^k-u^k)\, \sgn_n (\tilde{u}^k-u^k) \psi\\
&=
\int_{\Omega_r} \left| \nabla (\tilde{u}^k-u^k)\right|^2 \sgn_n' (\tilde{u}^k-u^k) \psi
+
\int_{\Omega_r} \sgn_n (\tilde{u}^k-u^k) \nabla (\tilde{u}^k-u^k)\cdot \nabla \psi \\
&\ge 
\int_{\Omega_r} \sgn_n (\tilde{u}^k-u^k) \nabla (\tilde{u}^k-u^k)\cdot \nabla \psi \\
&\xrightarrow[n\to \infty]{}
\int_{\Omega_r} \nabla \left| \tilde{u}^k-u^k \right| \cdot \nabla \psi
= 
-\int_{\Omega_r} \left| \tilde{u}^k-u^k \right| \Delta \psi
\ge 
-\max_{\Omega_r}{\left| \Delta \psi \right|} \int_{\Omega_r} \left| \tilde{u}^k-u^k \right| .
\end{align*}
Therefore, using a similar argument to that of Lemma~\ref{Lem_time_translation}, we obtain 
\eqref{tt_4} where $\bar{u}^k$ and $\bar{v}^k$ are replaced by $\tilde{u}^k$ and $\tilde{v}^k$, respectively, and $C$ is a constant independent of $k$ and $\xi$. 
Thus, we have 
\begin{align*}
&\ld \tilde{u}^k-u^k \rd_{L^1(\Omega_r)}(t)
+\ld \tilde{v}^k-v^k \rd_{L^1(\Omega_r)}(t) \\
& \quad \le
C \lr \ld \tilde{u}_0^k-u_0^k \rd_{L^1(\Omega_r')}
+\ld \tilde{v}_0^k-v_0^k \rd_{L^1(\Omega_r')} \rr. 
\end{align*}
Integrating this in $t$ yields the desired estimate for $s=1$. The estimate for $s\in (1,2)$ is obtained from \eqref{holder}. 
\endproof

\begin{lemma}
\label{Lem_RFK_cond2}
For all $s\in [1,2)$ and $\ep >0$, there exists $\omega \Subset Q_T$ such that  
\begin{align*}
\ld u \rd_{L^s(Q_T\setminus \omega)}+\ld v \rd_{L^s(Q_T\setminus \omega)}<\ep.
\end{align*}
\end{lemma}
\proof
It follows from H\"older's inequality that 
\begin{align*}
\int_{Q_T\setminus \omega} \left| u^k \right|^s 
\le 
&\lr
\int_{Q_T\setminus \omega} \left| u^k \right|^{s\frac{2}{s}}
\rr^\frac{s}{2}
\lr
\int_{Q_T\setminus \omega} 1^{\frac{2}{2-s}}
\rr^\frac{2-s}{2}\\
&\le 
\ld u^k \rd_{L^2(Q_T)}^s \mbox{meas}(Q_T\setminus \omega)^\frac{2-s}{2}. 
\end{align*}
Lemma~\ref{Lem_apriori_L2H1} implies the desired estimate. 
\endproof

\section{Proof of Theorem~\ref{Thm_conv} in $L^p$ ($p\in [1,2)$) space.} 
By virtue of	Lemmas~\ref{Lem_apriori_L2H1}, \ref{Lem_time_translation}, \ref{Lem_space_translation} and \ref{Lem_RFK_cond2}, and the Riesz-Fr\'echet-Kolmogorov theorem, there exist subsequences $\{u^{k_n}\}$ and $\{v^{k_n}\}$ of $\{u^k\}$ and $\{v^k\}$, respectively, and functions $u^*\in L^\infty(0,T;L^2(\Omega))$ $\cap L^2(0,T;H^1(\Omega))$ and $v^*\in L^\infty(0,T;L^2(\Omega))$ ($v^*\in  L^2(0,T;H^1(\Omega))$ if $d_2>0$) such that 
\begin{align}
&u^{k_n} \to u^* & & \mbox{strongly in } L^s(Q_T)\ \forall s\in [1,2), \mbox{ a.e. in } Q_T,   \notag \\
& & &  \mbox{weakly in } L^2(0,T;H^1(\Omega)), \label{conv_u}\\
&v^{k_n} \to v^* && \mbox{strongly in } L^s(Q_T)\ \forall s\in [1,2),  \mbox{ a.e. in } Q_T,  \notag \\
& & & \mbox{weakly in } L^2(Q_T)\ (\mbox{in } L^2(0,T;H^1(\Omega)) \mbox{ if } d_2>0) \label{conv_v}
\end{align}
as $k_n$ tends to infinity.
It follows from \eqref{conv_u} and \eqref{conv_v} that  
\begin{align}
f_i(u^{k_n},v^{k_n}) &\to f_i(u^{*},v^{*}) \quad (i=1,2),\label{conv_fi}\\
F(u^{k_n},v^{k_n}) &\to F(u^{*},v^{*}) \notag \\
 \alpha_0(u^{k_n})&\to \alpha_0(u^*)\notag
\end{align}
a.e. in $Q_T$. 
In view of \eqref{conv_fi} and $\ld f_i(u^{k_n},v^{k_n}) \rd_{L^2(Q_T)} \le C$, $f_i(u^{k_n},v^{k_n})$ converges weakly in $L^2(Q_T)$. 
Passing to the limit as $k_n\to \infty$ in \eqref{conv_alphav} (or in \eqref{conv_ugamma}), we obtain 
$F(u^*(x,t),v^*(x,t))=0$ or $\alpha_0(u^*(x,t))=v^*(x,t)$ for a.e. $(x,t)\in Q_T$. Both imply $v^*(x,t) \in \alpha (u^*(x,t))$. Therefore, we have 
\begin{equation}
\label{v-alpha}
v^* \in \alpha (u^*). 
\end{equation}
It follows from the weak form of the equations for $u^k$ and $v^k$ that 
\begin{align*}
-\int_0^T & \la \frac{\partial \zeta}{\partial t}, u^k +v^k \ra 
+\int_0^T\lr \nabla (d_1u^k+d_2v^k),\nabla \zeta \rr \\
&= 
\lr u_{0}^k +v_0^k,  \zeta(0) \rr
+
\int_0^T\lr (f_1+f_2)(u^k,v^k), \zeta \rr
\end{align*}
for all $\zeta \in H^1(0,T;H^1(\Omega)^*)\cap L^2(0,T;H^1(\Omega))$. 
Hence, passing to the limit along subsequences, we have 
\begin{align*}
-\int_0^T & \la \frac{\partial \zeta}{\partial t}, u^* +v^* \ra 
+\int_0^T\lr \nabla (d_1u^*+d_2v^*),\nabla \zeta \rr \\
&= 
\lr u_{0} +v_0,  \zeta(0) \rr
+
\int_0^T\lr (f_1+f_2)(u^*,v^*), \zeta \rr
\end{align*}
Set $z = u^*+v^*$, then $u^*=(I+\alpha)^{-1}z$, $v^*=(I-(I+\alpha)^{-1})z$ and $d_1u^*+d_2v^*=(d_2I+(d_1-d_2)(I+\alpha)^{-1})z$ by virtue of \eqref{v-alpha}. 
Therefore, $z$ is a weak solution of \eqref{limiteq}, that completes the proof in $L^s(Q_T)$ ($s\in [1,2)$). 
\endproof

\section{Boundedness and strong $L^p$ ($p\ge 1$) convergence.}

\begin{lemma}
\label{Lem_apriori_Lp}
Assume \eqref{assum_d}--\eqref{assum_F}, and $u_0,v_0\in L^p(\Omega)$ for $2\le p\le \infty$. Let $(u^k,v^k)$ be the weak solution of \eqref{rd} with an initial datum $(u_0,v_0)$. Then, there exists a positive constant $C$ independent of $k$ such that 
\begin{align*}
\| u^k \|_{L^\infty(0,T;L^p(\Omega))}
+
\| v^k \|_{L^\infty(0,T;L^p(\Omega))}
\le 
C. 
\end{align*}
\end{lemma}

\proof
The proof is similar to that of Lemma~\ref{Lem_apriori_L2H1}. 
Let $p \in [2,\infty)$ be fixed. We consider $u_0^{\ep},v_0^{\ep} \in L^p(\Omega)\cap H^1(\Omega)$ such that $u_0^{\ep} \to u_0$ and $v_0^{\ep} \to v_0$ in $L^p(\Omega)$ as $\ep$ tends to zero. 
Let $(u^{k,\ep},v^{k,\ep})$ be a  solution of \eqref{rd} with an initial datum $(u_0^{\ep},v_0^{\ep})$. 
Then, we have 
\[
{u}^{k,\ep} \to {u}^k ,\quad {v}^{k,\ep} \to {v}^k 
\quad \mbox{in } L^\infty(0,T;L^p(\Omega))
\]
as $\ep$ tends to zero.

Define functions $\varphi^p$ and $\Phi_g^p$ for $p \in [2,\infty)$ and for a given nondecreasing Lipschitz continuous function $g$ with a Lipschitz constant $L_g$ as follows. 
\[
\varphi^p(s)=|s|^{p-2}s, \qquad 
\Phi_g^p(s)=
\int_0^s \varphi^p(g(r))dr
\quad \mbox{for } s\in \R. 
\]
The following inequality can be obtained easily. 
\begin{equation*}
\frac{1}{pL_g} \left| g(s) \right|^p 
\le \Phi_g^p (s) \le 
\frac{L_g^{p-1}}{p} \left| s\right|^p \quad \mbox{for } s\in \R. 
\end{equation*}
We use the same functions $\gamma_\delta$ and $\alpha_\delta$ ($\delta \ge 0$) as those in the proof of Lemma~\ref{Lem_apriori_L2H1}. 
%
%
%
%
%
%
Considering 
\begin{align*}
\int_0^{t}\!\! \int_\Omega \Big\{ &
(\mbox{eq.}\, u^{k,\ep}) \times \varphi^p(u^{k,\ep})
+
(\mbox{eq.}\, v^{k,\ep}) \times \varphi^p\lr \gamma_\delta(v^{k,\ep})\rr \\
&+
(\mbox{eq.}\, u^{k,\ep}) \times \varphi^p \lr \alpha_\delta(u^{k,\ep}) \rr 
+
(\mbox{eq.}\, v^{k,\ep}) \times \varphi^p(v^{k,\ep}) 
\Big\}, 
\end{align*}
and following the same strategy as in the proof Lemma~\ref{Lem_apriori_L2H1}, we obtain 
\begin{align*}
 \ld u^{k,\ep} (t) \rd_{L^p(\Omega)}^p
&+
 \ld v^{k,\ep} (t) \rd_{L^p(\Omega)}^p \\
\le &
Ce^p \lr 
\ld u_0^{\ep} \rd_{L^p(\Omega)}^p
+
\ld v_0^{\ep} \rd_{L^p(\Omega)}^p
+
p\int_\Omega \Phi_{\alpha_0}^p(u_0^\ep)
+
p\int_\Omega \Phi_{\gamma_0}^p(v_0^\ep)
\rr. 
\end{align*}
Note that 
\begin{align*}
\Phi_{\alpha_0}^p (s)
&\le 
\frac{1}{C_3p}\left\{ \lr {C_3}|s|+{C_4} \rr^p - {C_4}^p \right\} \quad \mbox{for } s\in \R, \\
\Phi_{\gamma_0}^p (s)
&\le 
\frac{C_1}{p}\left\{ \lr \frac{1}{C_1}|s|+\frac{C_2}{C_1} \rr^p - \lr \frac{C_2}{C_1} \rr^p \right\} \quad \mbox{for } s\in \R, 
\end{align*}
where $C_i$ $(i=1,2,3,4)$ are constants given in \eqref{assum_F}. 
By virtue of the Minkowski inequality, the following estimate holds.
\begin{align*}
 \ld u^{k,\ep} (t) \rd_{L^p(\Omega)}
+
 \ld v^{k,\ep} (t) \rd_{L^p(\Omega)} 
\le 
CC^{1/p}p^{1/p} \lr 
1+ 
\ld u_0^{\ep} \rd_{L^p(\Omega)}
+
\ld v_0^{\ep} \rd_{L^p(\Omega)}
\rr.
\end{align*}
Here $C$ is a positive constant independent of $k$, $\ep$ and $p$. Thus, we can pass to the limit as $\ep \to 0$. Then, the desired estimate for $p\in [2,\infty)$ holds. 
The result for $p=\infty$ is obtained by taking to the limit as $p\to +\infty$. 
\endproof

\begin{corollary}
\label{Lpconv}
In addition to Assumptions~\eqref{assum_d}--\eqref{assum_init}, if \eqref{init_linf} holds, then the limit functions 
$u$ and $v$ belong to $L^\infty (Q_T)$ and 
subsequences of $\{u^k\}$ and $\{v^k\}$ converge strongly in $L^p(Q_T)$ for any $p\ge 1$. 
\end{corollary}
Thus, the remaining parts of Theorem~\ref{Thm_conv} have been proved.

\section{Concluding remarks}

In this paper, we have considered the fast reaction limit of \eqref{rd} under the general setting. 
Our setting is so general as to include a lot of problems in the literature. 
The key assumption is the existence of a maximal monotone $\alpha$ which satisfies \eqref{rel_F_alpha}. 
Thanks to this, $(I+\alpha)^{-1}$ is a non-decreasing contraction mapping, and thus, the limit equation~\eqref{limiteq} is well-defined. 
The maximal monotonicity assumption of $\alpha$ might be weakened. 
If $A=\{ (u,v) \ |\ v\in \alpha (u) \}$ and $B_c=\{ (u,v) \ |\ u+v=c \}$ intersect at one point for all $c\in \R$, the fast reaction limit might be considered. 
For example, if there exist a constant $\lambda <1$ and a multi-valued function $\alpha:\R \to 2^\R$ such that $\alpha$ satisfies \eqref{rel_F_alpha} and $\lambda I+\alpha$ is maximal monotone, then $(I+\alpha)^{-1}=((1-\lambda)I+(\lambda I+\alpha))^{-1}$ is still non-decreasing and Lipschitz continuous. 
Therefore, under this assumption, the limit equation~\eqref{limiteq} is well-defined and one can expect the convergence of the solution of the reaction-diffusion system~\eqref{rd}. 
However, the strategy in this paper can not be applied to this situation. 
The proof under this assumption remains an open problem. 
On the other hand, if $A$ and $B_c$ intersect at more than two points for some $c$, the limit equation~\eqref{limiteq} is a forward-backward equation, that is, ill-posed. 
In this case, Turing instability for \eqref{rd} should occur. The evidences have been obtained for specific choices of $F$ in \eqref{rd} \cite{jm,mps,mo}.



\medskip
Received xxxx 20xx; revised xxxx 20xx.
\medskip

\end{document}